\documentclass{article}

\title{Expressions for values of the gamma function}

\author{Raimundas Vid\=unas\footnote{Supported by the 21 Century COE
Programme "Development of Dynamic Mathematics with High Functionality" of
the Ministry of Education, Culture, Sports, Science and Technology of
Japan.}\\
 \em Kyushu University}

\newtheorem{theorem}{Theorem}[section]
\newtheorem{lemma}[theorem]{Lemma}

\usepackage{latexsym}

\newcommand{\hpg}[5]{{}_{#1}\mbox{\rm F}_{\!#2}\!
  \left(\left.{#3 \atop #4}\right| #5 \right) }

\newcommand{\proof}{{\bf Proof. }}
\newcommand{\qed}{\hfill $\Box$\\}
\newcommand{\equal}{&\!\!=\!\!&}
\newcommand{\equall}{\;\,=\;\,}
\newcommand{\dst}{\displaystyle}

\newcommand{\KK}[1]{\mbox{\bf K}\!\left(#1\right)}
\newcommand{\Exp}[1]{\exp\left(#1\right)}

\newcommand{\QQ}{\mbox{\bf Q}}

\newcommand{\ZZ}{\mbox{\bf Z}}
\newcommand{\GAMMA}[2]{\Gamma\!\left(\frac{#1}{#2}\right)}

\begin{document}

\maketitle

\begin{abstract} This paper presents expressions for gamma values at
rational points with the denominator dividing 24 or 60. These gamma values
are expressed in terms of 10 distinct gamma values and rational powers of
$\pi$ and a few real algebraic numbers. Our elementary list of formulas can
be conveniently used to evaluate, for example, algebraic Gauss
hypergeometric functions by the Gauss identity. Also, algebraic independence
of gamma values and their relation to the elliptic {\bf K} function are
briefly discussed.
\end{abstract}

\section{Introduction}

The gamma function \cite[Chapter~1]{specfaar} satisfies the difference
equation
\begin{equation} \label{differeq}
\Gamma(x+1) = x\,\Gamma(x),
\end{equation}
the Euler reflection formula
\begin{equation} \label{reflection}
\Gamma(x)\,\Gamma(1-x)=\frac{\pi}{\sin(\pi x)},
\end{equation}
and the Gauss multiplication formula
\begin{equation} \label{gmultiplic}
\Gamma(x)\,\Gamma\left(x+\frac1n\right)\ldots\Gamma\left(x+\frac{n\!-\!1}{n}\right)
= n^{\frac12-nx}\,(2\pi)^{\frac{n-1}2}\,\Gamma(nx).
\end{equation}
In the last formula, $n$ is a positive integer. Its special case $n=2$ is
known as Legendre's duplication formula. We refer to these functional
equations for the gamma function as the {\em standard equations}.

Values of the gamma function at rational points are of broad interest. By
historical motivation, $\Gamma(n)=(n-1)!$ when $n$ is a positive integer. An
easy consequence of the reflection formula is $\GAMMA12=\sqrt{\pi}$. By
using difference equation (\ref{differeq}) one can evaluate $\Gamma(x)$ for
rational $x$ with the denominator $2$. No explicit evaluations of other
gamma values are known. Some gamma terms (i.e., quotients of products of
gamma values) occur as values of hypergeometric functions at special points
\cite{specfaar} and as period integrals \cite{periods}, \cite{deligne}. In
particular, this applies to values of the elliptic $\bf K$-function at
so-called {\em elliptic integral singular values} \cite{zucker},
\cite{selbergc}, \cite{campbell}. Conversely, some gamma values at rational
points can be expressed in terms of elliptic integrals; see \cite{borwzu},
\cite{waldschmidt} and Section \ref{gammaelliptic} in this paper.

The purpose of this paper is to present explicit relations between gamma
values in the set
\begin{equation} \label{gammaset}
\left\{ \left. \GAMMA{k}{n} \;\right|\; k,n\in\ZZ;\ 0<\frac{k}{n}<1;\
\mbox{$n$ divides 24 or 60} \right\}.
\end{equation}
We show that standard formulas (\ref{reflection})--(\ref{gmultiplic}) imply
that all these gamma values can be multiplicatively expressed in terms of
rational powers of $\pi$, rational powers of few algebraic numbers, and the
following 10 gamma values:
\begin{equation} \label{gammabasis}
\begin{array}{ccccc}
\dst \GAMMA13, & \dst \GAMMA14,& \dst \GAMMA15,
& \dst \GAMMA25,& \dst \GAMMA18, \vspace{2pt}\\
\dst \GAMMA1{15},& \dst \GAMMA1{20}, & \dst \GAMMA1{24},& \dst \GAMMA1{60},
& \dst \GAMMA7{60}.
\end{array}
\end{equation}
In Section \ref{gammaevaluations} we present explicit expressions for the
gamma values in (\ref{gammaset}) in these terms. By using difference
equation (\ref{differeq}) and our list of formulas, one can express any
gamma value at a rational point with the denominator dividing 24 or 60 in
the same fashion. In Section \ref{generalresults} we indicate an elementary
proof of our evaluations, and show that the standard formulas do not imply
any relations between the distinguished values in (\ref{gammabasis}). If
Lang's conjecture \cite{lang} is true, those 10 gamma values and the
constant $\pi$ are algebraically independent over $\QQ$.

Relations between gamma values at rational points with small denominators
are widely known; see \cite[page~\sf GammaFunction.html]{wolfram} for
example. Some tricky relations between gamma values from the set
(\ref{gammaset}) are derived in \cite{borwzu}, \cite{kato} and probably by
other authors. Of special interest are gamma terms with algebraic values
\cite{koblitzo}, \cite{pdas}. Our list of relations between gamma values may
be useful for many purposes. For example, explicit evaluation of algebraic
Gauss hypergeometric functions at $x=1$ requires gamma values at rational
points with the denominators dividing 24 or 60; see Section
\ref{application}. Our formula list is precisely enough for this purpose. In
Section 5 we present available expressions of gamma values in
(\ref{gammabasis}) in terms of the elliptic $\bf K$-function.

\section{Explicit formulas}
\label{gammaevaluations}

Here we present a list of explicit expressions of gamma values in the set
(\ref{gammaset}) in terms of the gamma values in (\ref{gammabasis}).
A proof of these relations is indicated in the following Section. A
comparable list of expressions for gamma values is available at
\cite{webval}.

To make formulas more compact, we introduce the following constants:
\begin{equation} \label{cconstants}
\begin{array}{ll}
\phi=5+\sqrt5, & \phi^\star=5-\sqrt5,\vspace{2pt}\\
\psi=\sqrt{5+2\sqrt5},\qquad & \psi^\star=\sqrt{5-2\sqrt5}.
\end{array}
\end{equation}
Note that
\begin{equation}
\phi\,\phi^\star=20,\qquad \psi\,\psi^\star=\sqrt{5},\qquad
\psi=\frac{\phi\,^{3/2}}{2^{3/2}\,\sqrt{5}},\qquad
\psi^\star=\frac{\phi^\star\,^{3/2}}{2^{3/2}\,\sqrt{5}}.
\end{equation}
We express gamma values from (\ref{gammaset}) in terms of the values in
(\ref{gammabasis}) and rational powers of the following constants:
\begin{equation} \label{algnumbers}
\begin{array}{c}
\pi,\quad 2,\quad 3,\quad 5,\quad \sqrt2\pm1,\quad \sqrt3\pm 1,\quad
\sqrt{3}\pm\sqrt2,\quad \sqrt5\pm\sqrt3,\\
\phi,\quad \phi^\star,\quad\sqrt{15}\pm\psi,\quad
\sqrt{15}\pm\psi^\star,\quad \sqrt{10}\pm\sqrt{\phi},\quad
\sqrt{10}\pm\sqrt{\phi^\star}.
\end{array}
\end{equation}
When applying our formulas to gamma terms, it is easy to invert and simplify
these algebraic numbers. Some useful expressions with these algebraic
numbers 
are presented in Lemma \ref{algrelations} in the next Section.

Here is our list of formulas:
\begin{eqnarray*}
\GAMMA12\equal\sqrt{\pi}\,,\hspace{126pt}
\GAMMA23\equall\frac{2\,\pi}{\sqrt3}\;\GAMMA13^{-1},\\
\GAMMA34\equal\pi\,\sqrt2\,\;\GAMMA14^{-1},\hspace{76pt}
\GAMMA16\equall\frac{\sqrt3}{\sqrt\pi\;2^{1/3}}\;\GAMMA13^2,\\
\GAMMA35\equal\frac{\pi\,\sqrt2\,\sqrt{\phi^\star}}{\sqrt5}\;\GAMMA25^{-1},\hspace{54pt}
\GAMMA56\equall\frac{\pi^{3/2}\,2^{4/3}}{\sqrt3}\;\GAMMA13^{-2},\\
\GAMMA45\equal\frac{\pi\,\sqrt2\;\sqrt{\phi}}{\sqrt5}\;\GAMMA15^{-1},\hspace{58pt}
\GAMMA38\equall\sqrt{\pi}\,\sqrt{\sqrt2-1}\;\GAMMA14^{-1}\GAMMA18,\\
\GAMMA58\equal\sqrt{\pi}\,2^{3/4}\;\GAMMA14\GAMMA18^{-1},\hspace{35pt}
\GAMMA78\equall\pi\,2^{3/4}\,\sqrt{\sqrt2+1}\;\GAMMA18^{-1},\\
\GAMMA1{10}\equal\frac{\sqrt{\phi}}{\sqrt{\pi}\;2^{7/10}}\;\GAMMA15\GAMMA25,\hspace{34pt}
\GAMMA3{10}\equall\frac{\sqrt{\pi}\;\phi^\star}{2^{3/5}\,\sqrt{5}}\;\GAMMA15\GAMMA25^{-1},\\
\GAMMA7{10} \equal \sqrt{\pi}\;2^{3/5}\;\GAMMA15^{-1}\GAMMA25, \hspace{30pt}
\GAMMA9{10}\equall\frac{\pi^{3/2}\,2^{7/10}\,\sqrt{\phi}}{\sqrt{5}}\;\GAMMA15^{-1}\GAMMA25^{-1},\\
\GAMMA1{12} \equal \frac{3^{3/8}\,\sqrt{\sqrt3+1}}{\sqrt{\pi}\;2^{1/4}}\;\GAMMA13\GAMMA14,\\
\GAMMA5{12} \equal\frac{\sqrt{\pi}\;2^{1/4}\,\sqrt{\sqrt3-1}}{3^{1/8}}\;\GAMMA14\GAMMA13^{-1},\\
\GAMMA7{12} \equal\sqrt{\pi}\;2^{1/4}\,3^{1/8}\,\sqrt{\sqrt3-1}\;\GAMMA13\GAMMA14^{-1},\\
\GAMMA{11}{12} \equal\frac{\pi^{3/2}\,2^{3/4}\,\sqrt{\sqrt3+1}}{3^{3/8}}\;\GAMMA13^{-1}\GAMMA14^{-1},\\
\GAMMA2{15}\equal\frac{\sqrt{\phi^\star}\,\sqrt{\sqrt{15}-\psi^\star}}
{2\cdot 3^{7/20}\,5^{1/3}}\;\GAMMA13^{-1}\GAMMA25\GAMMA1{15},\\
\GAMMA4{15}\equal\frac{\sqrt{\phi}\,\sqrt{\sqrt{15}-\psi}\,\sqrt{\sqrt{15}-\psi^\star}}
{2^{3/2}\,3^{3/10}\,\sqrt{5}}\;\GAMMA15^{-1}\GAMMA25\GAMMA1{15},\\
\GAMMA7{15}\equal\frac{3^{9/20}\,\sqrt{\phi^\star}\,\sqrt{\sqrt{15}+\psi^\star}}
{2\cdot 5^{1/6}}\;\GAMMA13\GAMMA15\GAMMA1{15}^{-1},\\
\GAMMA8{15}\equal\frac{\pi\,\sqrt2\,\sqrt{\sqrt{15}-\psi}}{3^{9/20}\,5^{1/3}}
\;\GAMMA13^{-1}\GAMMA15^{-1}\GAMMA1{15},\\
\GAMMA{11}{15}\equal2\,\pi\cdot 3^{3/10}\;\GAMMA15\GAMMA25^{-1}\GAMMA1{15}^{-1},\\
\GAMMA{13}{15}\equal\frac{\pi\,\sqrt2\;3^{7/20}\,\sqrt{\sqrt{15}+\psi}}{5^{1/6}}
\;\GAMMA13\GAMMA25^{-1}\GAMMA1{15}^{-1},\\
\GAMMA{14}{15}\equal\frac{\pi\,\sqrt{\phi}\,\sqrt{\sqrt{15}+\psi}\,\sqrt{\sqrt{15}+\psi^\star}}
{\sqrt{2}\,\sqrt{5}}\;\GAMMA1{15}^{-1},\\
\GAMMA3{20}\equal\frac{\sqrt\pi\;\phi^\star\,\sqrt{\sqrt{10}-\sqrt{\phi^\star}}}
{2^{21/20}\,5^{7/8}}\;\GAMMA25^{-1}\GAMMA1{20},\\
\GAMMA7{20}\equal\frac{\sqrt\pi\;\sqrt{\sqrt{10}-\sqrt{\phi}}}
{2^{3/20}\,5^{3/8}}\;\GAMMA15^{-1}\GAMMA1{20},\\
\GAMMA9{20}\equal\frac{\pi\sqrt{\sqrt{10}-\sqrt{\phi}}\;\sqrt{\sqrt{10}-\sqrt{\phi^\star}}}
{2^{1/5}\,\sqrt5}\;\GAMMA15^{-1}\GAMMA25^{-1}\GAMMA1{20},\\
\GAMMA{11}{20}\equal2^{1/5}\,\sqrt\phi\;\GAMMA15\GAMMA25\GAMMA1{20}^{-1},\\
\GAMMA{13}{20}\equal\frac{\sqrt\pi\;2^{3/20}\,\sqrt{\phi^\star}\,
\sqrt{\sqrt{10}+\sqrt{\phi^\star}}}{5^{1/8}}\;\GAMMA15\GAMMA1{20}^{-1},\\
\GAMMA{17}{20}\equal\frac{\sqrt\pi\;2^{1/20}\,\sqrt{\phi}\,\sqrt{\sqrt{10}+\sqrt{\phi}}}
{5^{1/8}}\;\GAMMA25\GAMMA1{20}^{-1},\\
\GAMMA{19}{20}\equal\frac{\pi\,\sqrt\phi\,\sqrt{\sqrt{10}+\sqrt{\phi}}\;\sqrt{\sqrt{10}+\sqrt{\phi^\star}}}
{\sqrt5}\;\GAMMA1{20}^{-1},\\
\GAMMA5{24}\equal\frac{\sqrt\pi\;\sqrt{\sqrt2-1}\,\sqrt{\sqrt3-1}}{2^{1/6}\,\sqrt3}
\;\GAMMA13^{-1}\GAMMA1{24},\\
\GAMMA7{24}\equal\frac{\sqrt\pi\,\sqrt{\sqrt3-1}\,\sqrt{\sqrt3-\sqrt2}}
{2^{1/4}\,3^{3/8}}\;\GAMMA14^{-1}\GAMMA1{24},\\
\GAMMA{11}{24}\equal\frac{\pi\,2^{1/12}\,\sqrt{\sqrt2-1}\,\sqrt{\sqrt3-\sqrt2}}
{3^{3/8}}\;\GAMMA13^{-1}\GAMMA14^{-1}\GAMMA1{24},\\
\GAMMA{13}{24}\equal 2^{2/3}\,3^{3/8}\,\sqrt{\sqrt3+1}\;\GAMMA13\GAMMA14\GAMMA1{24}^{-1},\\
\GAMMA{17}{24}\equal 2\,\sqrt\pi\;3^{3/8}\,\sqrt{\sqrt2+1}\;\GAMMA14\GAMMA1{24}^{-1},\\
\GAMMA{19}{24}\equal\sqrt\pi\;2^{11/12}\,\sqrt3\,\sqrt{\sqrt3+\sqrt2}\;\GAMMA13\GAMMA1{24}^{-1},\\
\GAMMA{23}{24}\equal\pi\;2^{3/4}\,\sqrt{\sqrt2+1}\,\sqrt{\sqrt3+1}\,\sqrt{\sqrt3+\sqrt2}\;\GAMMA1{24}^{-1},\\
\GAMMA1{30}\equal\frac{3^{9/20}\,\sqrt{\phi}\,\sqrt{\sqrt{15}+\psi}}
{\sqrt\pi\;2^{16/15}\,5^{1/6}}\;\GAMMA13\GAMMA15,\\
\GAMMA7{30}\equal\frac{3^{3/20}\,\sqrt{\phi^\star}\,\sqrt{\sqrt{15}+\psi^\star}}
{\sqrt\pi\;2^{22/15}\,5^{1/6}}\;\GAMMA13\GAMMA25,\\
\GAMMA{11}{30}\equal\frac{\sqrt\pi\,\sqrt{\phi}\,\sqrt{\sqrt{15}-\psi}}
{2^{11/15}\,3^{1/20}\,5^{1/3}}\;\GAMMA13^{-1}\GAMMA15,\\
\GAMMA{13}{30}\equal\frac{\sqrt\pi\;3^{7/20}\,\phi^\star\,\sqrt{\sqrt{15}-\psi^\star}}
{2^{41/30}\,5^{2/3}}\;\GAMMA13\GAMMA25^{-1},\\
\GAMMA{17}{30}\equal\frac{\sqrt\pi\;\sqrt{\phi^\star}\,\sqrt{\sqrt{15}-\psi^\star}}
{2^{2/15}\,3^{7/20}\,5^{1/3}}\;\GAMMA13^{-1}\GAMMA25,\\
\GAMMA{19}{30}\equal\frac{\sqrt\pi\;3^{1/20}\,\phi\,\sqrt{\sqrt{15}-\psi}}
{2^{23/30}\,5^{2/3}}\;\GAMMA13\GAMMA15^{-1},\\
\GAMMA{23}{30}\equal\frac{\pi^{3/2}\,\phi^\star\,\sqrt{\sqrt{15}+\psi^\star}}
{2^{1/30}\,3^{3/20}\,5^{5/6}}\;\GAMMA13^{-1}\GAMMA25^{-1},\\
\GAMMA{29}{30}\equal\frac{\pi^{3/2}\,\phi\,\sqrt{\sqrt{15}+\psi}}
{2^{13/30}\,3^{9/20}\,5^{5/6}}\;\GAMMA13^{-1}\GAMMA15^{-1},\\
\GAMMA{11}{60}\equal\frac{\sqrt\pi\;\sqrt{\phi}\;\sqrt{\sqrt{15}-\psi}\,\sqrt{\sqrt{10}-\sqrt{\phi}}}
{2^{5/4}\,\sqrt3\;5^{17/24}}\;\GAMMA13^{-1}\GAMMA1{60},\\
\GAMMA{13}{60}\equal\frac{\sqrt\pi\;\sqrt{\phi^\star}\,\sqrt{\sqrt3+1}\,\sqrt{\sqrt5-\sqrt3}
\,\sqrt{\sqrt{15}-\psi^\star}}{2^{13/10}\,3^{3/20}\,5^{3/8}}\;\GAMMA25^{-1}\GAMMA7{60},\\
\GAMMA{17}{60}\equal\frac{\sqrt\pi\;\sqrt{\phi^\star}\,\sqrt{\sqrt{15}-\psi^\star}\,
\sqrt{\sqrt{10}-\sqrt{\phi^\star}}}{2^{3/4}\,\sqrt3\;5^{11/24}}\;\GAMMA13^{-1}\GAMMA7{60},\\
\GAMMA{19}{60}\equal\frac{\sqrt\pi\;\sqrt{\phi}\;\sqrt{\sqrt3-1}\,\sqrt{\sqrt5-\sqrt3}\,
\sqrt{\sqrt{15}-\psi}}{2^{7/5}\,3^{9/20}\,5^{5/8}}\;\GAMMA15^{-1}\GAMMA1{60},\\
\GAMMA{23}{60}\equal\frac{\pi\,\sqrt{\phi^\star}\,\sqrt{\sqrt3+1}\,\sqrt{\sqrt5-\sqrt3}\,
\sqrt{\sqrt{10}-\sqrt{\phi^\star}}}{2^{11/20}\,3^{3/20}\,5^{7/12}}\;\GAMMA13^{-1}\GAMMA25^{-1}\GAMMA7{60},\\
\GAMMA{29}{60}\equal\frac{\pi\,\sqrt{\phi}\;\sqrt{\sqrt3-1}\,\sqrt{\sqrt5-\sqrt3}\,
\sqrt{\sqrt{10}-\sqrt{\phi}}}{2^{23/20}\,3^{9/20}\,5^{7/12}}\;\GAMMA13^{-1}\GAMMA15^{-1}\GAMMA1{60},\\
\GAMMA{31}{60}\equal\frac{3^{9/20}\,\sqrt\phi\;\sqrt{\sqrt{15}+\psi}}{2^{1/10}\,5^{1/6}}
\;\GAMMA13\GAMMA15\GAMMA1{60}^{-1},\\
\GAMMA{37}{60}\equal\frac{3^{3/20}\,\sqrt{\phi^\star}\,\sqrt{\sqrt{15}+\psi^\star}}
{2^{7/10}\,5^{1/6}}\;\GAMMA13\GAMMA25\GAMMA7{60}^{-1},\\
\GAMMA{41}{60}\equal\frac{\sqrt\pi\;2^{3/20}\,3^{9/20}\,\sqrt\phi\;\sqrt{\sqrt{10}+\sqrt\phi}}
{5^{1/8}}\;\GAMMA15^{-1}\GAMMA1{60},\\
\GAMMA{43}{60}\equal\frac{\sqrt\pi\,\sqrt3\;\sqrt{\phi^\star}\;\sqrt{\sqrt3-1}\,
\sqrt{\sqrt5+\sqrt3}}{\sqrt2\; 5^{7/24}}\;\GAMMA13\GAMMA7{60}^{-1},\\
\GAMMA{47}{60}\equal\frac{\sqrt\pi\;2^{1/20}\,3^{3/20}\,\sqrt{\phi^\star}\;
\sqrt{\sqrt{10}+\sqrt{\phi^\star}}}{5^{3/8}}\;\GAMMA25^{-1}\GAMMA7{60},\\
\GAMMA{49}{60}\equal\frac{\sqrt\pi\,\sqrt3\;\sqrt\phi\;\sqrt{\sqrt3+1}\,
\sqrt{\sqrt5+\sqrt3}}{5^{1/24}}\;\GAMMA13\GAMMA1{60}^{-1},\\
\GAMMA{53}{60}\equal\frac{\pi\,\phi^\star\sqrt{\sqrt3-1}\,\sqrt{\sqrt{5}+\sqrt{3}}\,
\sqrt{\sqrt{15}+\psi^\star}\,\sqrt{\sqrt{10}+\sqrt{\phi^\star}}}{2^{5/4}\,5^{3/4}}\;\GAMMA7{60}^{-1},\\
\GAMMA{59}{60}\equal\frac{\pi\,\phi\;\sqrt{\sqrt3+1}\,\sqrt{\sqrt{5}+\sqrt{3}}\;
\sqrt{\sqrt{15}+\psi}\;\sqrt{\sqrt{10}+\sqrt\phi}}{2^{5/4}\,5^{3/4}}\;\GAMMA1{60}^{-1}.
\end{eqnarray*}

\section{Proof of the formulas}
\label{generalresults}

Here we give an elementary proof of all identities in the previous Section.
We do not reproduce detailed computations in our proofs; they are quite
straightforward though tedious. Lemma \ref{algrelations} gives compact
expressions for some products of employed algebraic numbers
(\ref{algnumbers}). We also show that the standard equations do not imply
any relations between gamma values in (\ref{gammabasis}), and discuss
briefly algebraic independence of those gamma values. Along the way, we give
a set of 16 gamma values which are enough to express any gamma value at a
rational point with the denominator dividing 120.

First we indicate relevant values of the sine function. To read off the
sinus values, one has to compare formulas in Lemma \ref{cosecants} with the
following form of Euler's identity:
\begin{equation}
\Exp{ix}= \sin\left(\frac\pi2-x\right)+i\sin x. 
\end{equation}
For completeness, recall the well-known values
$\Exp{\frac{i\pi}{3}}=\frac{1+i\sqrt3}2$ and
$\Exp{\frac{i\pi}{4}}=\frac{1+i}{\sqrt2}$.
\begin{lemma} \label{cosecants}
With the same notation as in $(\ref{cconstants})$, the following formulas
hold:
\begin{eqnarray*}
\Exp{\frac{i\pi}{5}}\equal\phi\;\frac{1+i\,\psi^\star}{4\,\sqrt5},\hspace{78pt}
\Exp{\frac{2i\pi}{5}}\equall\phi^\star\,\frac{1+i\,\psi}{4\,\sqrt5},\\
\Exp{\frac{i\pi}{8}}\equal\frac{\sqrt{\sqrt2+1}+i\sqrt{\sqrt2-1}}{2^{3/4}},\hspace{28pt}
\Exp{\frac{i\pi}{12}}\equall\frac{\sqrt3+1+i\left(\sqrt3-1\right)}{2\,\sqrt2},\\
\Exp{\frac{i\pi}{15}}\equal
\phi^\star\;\frac{\sqrt3\,\psi+1+i\left(\psi-\sqrt3\right)}{8\,\sqrt5},\quad
\Exp{\frac{2i\pi}{15}}\equall
\phi\;\frac{\sqrt3\,\psi^\star+1+i\left(\psi^\star-\sqrt3\right)}{8\,\sqrt5},\\
\Exp{\frac{4i\pi}{15}}\equal
\phi^\star\;\frac{\sqrt3\,\psi-1+i\left(\psi+\sqrt3\right)}{8\,\sqrt5},\quad
\Exp{\frac{7i\pi}{15}}\equall
\phi\;\frac{\sqrt3\,\psi^\star-1+i\left(\sqrt3-\psi^\star\right)}{8\,\sqrt5},\\
\Exp{\frac{i\pi}{20}}\equal
\sqrt{\phi^\star}\;\frac{\psi+\sqrt5+i\left(\psi-\sqrt5\right)}{4\,\sqrt5},\quad
\Exp{\frac{3i\pi}{20}}\equall
\sqrt{\phi}\;\frac{\sqrt5+\psi^\star+i\left(\sqrt5-\psi^\star\right)}{4\,\sqrt5},\\
\Exp{\frac{i\pi}{24}}\equal\frac{\sqrt{2\sqrt2+\sqrt3+1}+i\,\sqrt{2\sqrt2-\sqrt3-1}}{2^{5/4}},\\
\Exp{\frac{7i\pi}{24}}\equal\frac{\sqrt{2\sqrt2-\sqrt3+1}+i\,\sqrt{2\sqrt2+\sqrt3-1}}{2^{5/4}},\\
\Exp{\frac{i\pi}{60}}\equal
\phi^\star\;\frac{\left(\sqrt3+1\right)\left(\psi-\sqrt3+2\right)+i
\left(\sqrt3-1\right)\left(\sqrt3-\psi+2\right)}{8\,\sqrt{10}},\\
\Exp{\frac{7i\pi}{60}}\equal
\phi\;\frac{\left(\sqrt3-1\right)\left(\psi^\star\!+\sqrt3+2\right)+i
\left(\sqrt3+1\right)\left(\sqrt3+\psi^\star\!-2\right)}{8\,\sqrt{10}},\\
\Exp{\frac{13i\pi}{60}}\equal
\phi\;\frac{\left(\sqrt3+1\right)\left(\psi^\star\!-\sqrt3+2\right)+i
\left(\sqrt3-1\right)\left(\sqrt3-\psi^\star\!+2\right)}{8\,\sqrt{10}},\\
\Exp{\frac{19i\pi}{60}}\equal
\phi^\star\;\frac{\left(\sqrt3-1\right)\left(\psi+\sqrt3+2\right)+i
\left(\sqrt3+1\right)\left(\sqrt3+\psi-2\right)}{8\,\sqrt{10}}.\\
\end{eqnarray*}
\end{lemma}
\proof One can deduce the first two formulas by checking the sign of real
and imaginary parts on their right-hand sides, and checking the respective
equations $x^5\pm 1=0$. The remaining formulas can be consequently deduced
by inspecting the values of
\begin{eqnarray*}
\Exp{\frac{i\pi}{8}}^2,\ \Exp{\frac{i\pi}{12}}\Exp{\frac{i\pi}{4}},\
\Exp{\frac{i\pi}{15}}\Exp{\frac{i\pi}{3}},\ \Exp{\frac{i\pi}{15}}^2,\
\Exp{\frac{2i\pi}{15}}^2,\\ \Exp{\frac{2i\pi}{15}}\Exp{\frac{i\pi}{3}},\quad
\Exp{\frac{i\pi}{20}}\Exp{\frac{i\pi}{5}},\quad
\Exp{\frac{3i\pi}{20}}\Exp{\frac{i\pi}{4}},\\
\Exp{\frac{i\pi}{24}}^2,\quad\Exp{\frac{i\pi}{24}}\Exp{\frac{i\pi}{4}},\quad
\Exp{\frac{i\pi}{60}}\Exp{\frac{i\pi}4},\\
\Exp{\frac{i\pi}{12}}\Exp{\frac{i\pi}{5}},
\quad\Exp{\frac{i\pi}{60}}\Exp{\frac{i\pi}5},\quad\Exp{\frac{7i\pi}{60}}\Exp{\frac{i\pi}5}.
\end{eqnarray*}
Simplification of each exponential expression relates the respective
exponent value to earlier concluded exponent values, so it is enough to
check those relations with substituted algebraic expressions (and sometimes
the sign of the new exponent value).\qed

Now we present some expressions with employed algebraic numbers. They can be
used to simplify evaluated gamma terms and to transform the sinus values in
the previous Lemma.
\begin{lemma} \label{algrelations}
With the same notation as in $(\ref{cconstants})$, the following formulas
hold:
\begin{eqnarray*}
\sqrt{\sqrt{15}+\psi}\,\sqrt{\sqrt{15}-\psi}\equal\sqrt2\;\sqrt{\phi^\star},\hspace{46pt}
\sqrt{\sqrt{10}+\sqrt{\phi}}\,\sqrt{\sqrt{10}-\sqrt{\phi}}\equall\sqrt{\phi^\star},\\
\sqrt{\sqrt{15}+\psi^\star}\,\sqrt{\sqrt{15}-\psi^\star}\equal\sqrt2\;\sqrt{\phi},\hspace{42pt}
\sqrt{\sqrt{10}+\sqrt{\phi^\star}}\,\sqrt{\sqrt{10}-\sqrt{\phi^\star}}\equall\sqrt{\phi},\\
\sqrt{\sqrt{15}\pm\psi}\,\sqrt{\sqrt{15}\pm\psi^\star}\equal\frac{\sqrt{\phi^\star}}{\sqrt{2}}
\left(\psi\pm\sqrt3\right),\quad
\sqrt{\sqrt{10}\pm\sqrt{\phi}}\,\sqrt{\sqrt{10}\pm\sqrt{\phi^\star}}\equall\psi\pm\sqrt5,\\
\sqrt{\sqrt{15}\pm\psi}\,\sqrt{\sqrt{15}\mp\psi^\star}\equal\frac{\sqrt{\phi}}{\sqrt{2}}
\left(\sqrt3\pm\psi^\star\right),\quad
\sqrt{\sqrt{10}\pm\sqrt{\phi}}\,\sqrt{\sqrt{10}\mp\sqrt{\phi^\star}}\equall\sqrt5\pm\psi^\star,\\
\frac{\psi}{\sqrt5}\left(\sqrt{15}\pm\psi^\star\right)\equal\sqrt{3}\;\psi\pm1,\hspace{92pt}
\frac{\psi^\star}{\sqrt5}\left(\sqrt{15}\pm\psi\right)\equall\sqrt{3}\;\psi^\star\pm1,\\
\end{eqnarray*}\vspace{-30pt}
\begin{eqnarray*}
\hspace{10pt}
\sqrt{\sqrt{15}\pm\psi}\,\sqrt{\sqrt{10}\pm\sqrt{\phi}}\equal
\frac{\phi^\star}{2^{7/4}\,5^{1/4}}\,\sqrt{\sqrt3-1}\,\sqrt{\sqrt5+\sqrt3}\,\left(\sqrt3+2\pm\psi\right),\\
\sqrt{\sqrt{15}\mp\psi}\,\sqrt{\sqrt{10}\pm\sqrt{\phi}}\equal\frac{\phi^\star}{2^{7/4}\,5^{1/4}}\,
\sqrt{\sqrt3+1}\,\sqrt{\sqrt5-\sqrt3}\,\left(\psi\mp\sqrt3\pm2\right),\\
\sqrt{\sqrt{15}\pm\psi^\star}\,\sqrt{\sqrt{10}\pm\sqrt{\phi^\star}}\equal\frac{\phi}{2^{7/4}\,5^{1/4}}\,
\sqrt{\sqrt3+1}\,\sqrt{\sqrt5+\sqrt3}\,\left(\psi^\star\mp\sqrt3\pm 2\right),\\
\sqrt{\sqrt{15}\mp\psi^\star}\,\sqrt{\sqrt{10}\pm\sqrt{\phi^\star}}\equal
\frac{\phi}{2^{7/4}\,5^{1/4}}\,\sqrt{\sqrt3-1}\,\sqrt{\sqrt5-\sqrt3}\,\left(\sqrt3+2\pm\psi^\star\right),\\
2\sqrt2\pm\sqrt3+1\equal\left(\sqrt2+1\right)\left(\sqrt3\mp1\right)\left(\sqrt3\pm\sqrt2\right),\\
2\sqrt2\pm\sqrt3-1\equal\left(\sqrt2-1\right)\left(\sqrt3\pm1\right)\left(\sqrt3\mp\sqrt2\right),\\
2\sqrt3\pm\sqrt5\pm1\equal\frac{\phi^\star}{2\,\sqrt{5}}\;\left(\sqrt3\pm1\right)\left(\sqrt5\pm\sqrt3\right),\\
2\sqrt3\pm\sqrt5\mp1\equal\frac{\phi}{2\,\sqrt{5}}\;\left(\sqrt3\mp1\right)\left(\sqrt5\pm\sqrt3\right).
\end{eqnarray*}
Formulas with $\pm$, $\mp$ signs represent two identities, which can be read
by taking the upper signs or the lower signs respectively.
\end{lemma}
\proof The identities can be proved by direct manipulation. \qed

The following Lemma is important for proving independence of the gamma
values in (\ref{gammabasis}) with respect to the standard equations.
\begin{lemma} \label{kubert}
Let $N$ denote an integer greater than 2, and let $\varphi(N)$ denote
Euler's totient value. Let $\Sigma$ be the set of gamma values at rational
points with the denominator dividing $N$. Suppose that $\Sigma_0$ is a
minimal subset of $\Sigma$ which determines all other gamma values in
$\Sigma$ by the standard equations. Then $\Sigma_0$ has precisely
$\varphi(N)/2$ elements.
\end{lemma}
\proof This is a reformulation of Kubert's theorem in \cite{kubert}. See
also \cite[Chapter 2]{lang}. \qed

In the setting of the last Lemma, the set
$\left\{\left.\GAMMA{k}{N}\;\right|\; \gcd(k,N)=1,\ k<\frac{N}2\right\}$
looks like a natural candidate for $\Sigma_0$. This is a false impression in
general. One can check the formulas in Section \ref{gammaevaluations} and
notice that for $N=20,24,30,60$ gamma values at rational points with the
denominator $N$ depend on fewer gamma values of (\ref{gammabasis}) than
$\varphi(N)/2$. In these cases, gamma values in the mentioned set are
dependent.

Here we prove our main results. Along the way, we complement the set in
(\ref{gammabasis}) to a generating set for gamma values at rational points
with the denominator dividing 120. Note that $120$ is the lowest common
multiple of 24 and 60.
\begin{theorem} \label{mainresults}
\begin{enumerate}
\item[(1.)] Formulas in Section $\ref{gammaevaluations}$ hold.
 \item[(2.)] Standard equations $(\ref{differeq})$--$(\ref{gmultiplic})$ imply that
any gamma value with the denominator dividing $120$ can be expressed in
terms of $(\ref{gammabasis})$ and the following $6$ gamma values:
\begin{equation} \label{extravals}
\GAMMA1{40},\quad \GAMMA3{40},\quad \GAMMA7{40},\quad \GAMMA1{120},\quad
\GAMMA7{120},\quad \GAMMA{11}{120}.
\end{equation}
 \item[(3.)] Standard equations $(\ref{differeq})$--$(\ref{gmultiplic})$ do not
imply any relations\footnote{This statement also holds for the total set of
16 gamma values in (\ref{gammabasis}) and (\ref{extravals}).} between gamma
values in $(\ref{gammabasis})$.
\end{enumerate}
\end{theorem}
\proof For $x\in\QQ$, let $R(x)$ denote reflection formula
(\ref{reflection}), 
and let $M_n(x)$ denote multiplication formula (\ref{gmultiplic}). Consider
the following sequence of formulas:
\begin{eqnarray*} \label{eqs1}
R\!\left(\frac12\right), R\!\left(\frac13\right), R\!\left(\frac14\right),
R\!\left(\frac15\right), R\!\left(\frac25\right), M_2\!\left(\frac16\right),
R\!\left(\frac16\right), M_2\!\left(\frac18\right), R\!\left(\frac18\right),
R\!\left(\frac38\right),\\
M_2\!\left(\frac15\right),\ M_2\!\left(\frac1{10}\right),\
R\!\left(\frac1{10}\right),\ R\!\left(\frac3{10}\right),\
M_3\!\left(\frac1{12}\right),\ M_2\!\left(\frac1{12}\right),\
R\!\left(\frac1{12}\right),\ R\!\left(\frac5{12}\right),\\
M_3\!\left(\frac1{15}\right),\ R\!\left(\frac1{15}\right),\
R\!\left(\frac4{15}\right),\ M_5\!\left(\frac1{15}\right),\
M_3\!\left(\frac2{15}\right),\ R\!\left(\frac2{15}\right),\
R\!\left(\frac7{15}\right),\\
M_2\!\left(\frac1{20}\right),\ R\!\left(\frac1{20}\right),\
R\!\left(\frac9{20}\right),\ M_5\!\left(\frac1{20}\right),\
M_2\!\left(\frac3{20}\right),\ R\!\left(\frac3{20}\right),\
R\!\left(\frac7{20}\right),\\
M_3\!\left(\frac1{24}\right),\ M_2\!\left(\frac1{24}\right),\
M_2\!\left(\frac5{24}\right),\ R\!\left(\frac1{24}\right),\
R\!\left(\frac5{24}\right),\ R\!\left(\frac7{24}\right),\
R\!\left(\frac{11}{24}\right),\\
M_2\!\left(\frac1{15}\right),\ M_2\!\left(\frac2{15}\right),\
M_2\!\left(\frac1{30}\right),\ M_2\!\left(\frac7{30}\right),\
R\!\left(\frac1{30}\right),\ R\!\left(\frac7{30}\right),\
R\!\left(\frac{11}{30}\right),\ R\!\left(\frac{13}{30}\right),\\
M_3\!\left(\frac1{60}\right),\ M_3\!\left(\frac7{60}\right),\
M_2\!\left(\frac1{60}\right),\ M_2\!\left(\frac7{60}\right),\
M_2\!\left(\frac{11}{60}\right),\ R\!\left(\frac{13}{60}\right),\
M_2\!\left(\frac{13}{60}\right),\\
R\!\left(\frac1{60}\right),\ R\!\left(\frac7{60}\right),\
R\!\left(\frac{11}{60}\right),\ R\!\left(\frac{17}{60}\right),\
R\!\left(\frac{19}{60}\right),\ R\!\left(\frac{23}{60}\right),\
R\!\left(\frac{29}{60}\right).
\end{eqnarray*}
We claim that these formulas, seen as equations in gamma values from the set
(\ref{gammaset}), have a unique solution if the gamma values in
(\ref{gammabasis}) are assumed to be known. To prove this claim, note that
the first 14 formulas, up till $R\!\left(\frac3{10}\right)$, determine
consequently the gamma values at
\begin{equation} \label{gamvals1}
\frac12,\ \frac23,\ \frac34,\ \frac45,\ \frac35,\ \frac16,\ \frac56,\
\frac58,\ \frac78,\ \frac38,\ \frac7{10},\ \frac1{10},\ \frac9{10},\
\frac3{10}.
\end{equation}
Each of those 14 formulas relates the respective gamma value to earlier
concluded gamma values and the gamma values in (\ref{gammabasis}), like in
the proof of Lemma \ref{cosecants}. The next 4 formulas, up till
$R\!\left(\frac5{12}\right)$, express the products $\GAMMA1{12}\GAMMA5{12}$,
$\GAMMA1{12}\GAMMA7{12}$, $\GAMMA1{12}\GAMMA{11}{12}$,
$\GAMMA5{12}\GAMMA{7}{12}$  in terms of gamma values in (\ref{gammabasis})
and at the points in (\ref{gamvals1}). Straightforward manipulation of those
4 formulas expresses the gamma values at $\frac1{12}, \frac5{12},
\frac7{12}, \frac{11}{12}$ in terms of (\ref{gammabasis}) and the earlier
values. Similarly, the next 3 formulas, up till
$R\!\left(\frac4{15}\right)$, consequently determine $\GAMMA{11}{15}$,
$\GAMMA{14}{15}$, $\GAMMA4{15}$; and then the subsequent 4 formulas all
together determine $\GAMMA2{15}$, $\GAMMA7{15}$, $\GAMMA8{15}$,
$\GAMMA{13}{15}$. In the same way, the next formulas
$M_2\!\left(\frac1{20}\right)$, $R\!\left(\frac1{20}\right)$,\
$R\!\left(\frac9{20}\right)$ consequently determine $\GAMMA{11}{20}$,
$\GAMMA{19}{20}$, $\GAMMA9{20}$; and then the following 4
formulas\footnote{Gamma values cannot be determined entirely consequently.
As was noted in \cite{pdas}, some relations between gamma values are implied
by the standard equations up to taking square roots of both sides of an
equality.} all together determine $\GAMMA3{20}$, $\GAMMA7{20}$,
$\GAMMA{13}{20}$, $\GAMMA{17}{20}$. The remaining formulas, starting from
$M_3\!\left(\frac1{24}\right)$, consequently determine gamma values at
\begin{eqnarray*}
&\dst \frac{17}{24},\ \frac{13}{24},\ \frac5{24},\ \frac{23}{24},\
\frac{19}{24},\ \frac7{24},\ \frac{11}{24},\ \frac{17}{30},\ \frac{19}{30},\
\frac1{30},\ \frac7{30},\ \frac{29}{30},\ \frac{23}{30},\ \frac{11}{30},\ \frac{13}{30},\\
&\dst \frac{41}{60},\ \frac{47}{60},\ \frac{31}{60},\ \frac{37}{60},\
\frac{11}{60},\ \frac{13}{60},\ \frac{43}{60},\ \frac{59}{60},\
\frac{53}{60},\ \frac{49}{60},\ \frac{17}{60},\ \frac{19}{60},\
\frac{23}{60},\ \frac{29}{60}.
\end{eqnarray*}
We see that formulas are independent, and that they unique determine the
mentioned gamma values. Also, the gamma values in (\ref{gammabasis}) and the
consequently concluded gamma values exhaust all elements of the set in
(\ref{gammaset}). The claim follows.

To prove the first statement of the theorem, it is enough now to check that
evaluations in Section \ref{gammaevaluations} are compatible with the
sequence of formulas introduced at the beginning of this proof. Identities
in Lemmas \ref{cosecants} and \ref{algrelations} are very helpful in these
computations.

Here we show the second statement of the Theorem. We need to check the gamma
values at rational points with the denominators 40 and 120. The equations
\[
M_2\!\left(\frac1{40}\right),\ M_2\!\left(\frac3{40}\right),\
M_2\!\left(\frac7{40}\right),\ R\!\left(\frac{19}{40}\right),\
M_5\!\left(\frac3{40}\right),\ M_2\!\left(\frac{11}{40}\right)
\]
determine consequently the gamma values at $\frac{21}{40}$, $\frac{23}{40}$,
$\frac{27}{40}$, $\frac{19}{40}$, $\frac{11}{40}$, $\frac{31}{40}$ in terms
of (\ref{gammabasis}) and (\ref{extravals}). After this, other gamma values
at rational points with the denominator 40 are determined by reflection
formulas. Similarly, the equations
\begin{eqnarray*}
M_3\!\left(\frac1{120}\right),\ M_3\!\left(\frac7{120}\right),\
M_3\!\left(\frac{11}{120}\right),\ M_2\!\left(\frac1{120}\right),\
M_2\!\left(\frac7{120}\right),\ M_2\!\left(\frac{11}{120}\right),\
M_2\!\left(\frac{31}{120}\right),\\
M_2\!\left(\frac{41}{120}\right),\ M_2\!\left(\frac{47}{120}\right),\
R\!\left(\frac{59}{120}\right),\ M_5\!\left(\frac{11}{120}\right),\
M_3\!\left(\frac1{40}\right),\ M_2\!\left(\frac{23}{120}\right),\
M_2\!\left(\frac{43}{120}\right)
\end{eqnarray*}
determine consequently the gamma values at
\begin{equation} \label{gamvals2}
\frac{41}{120},\ \frac{47}{120},\ \frac{91}{120},\ \frac{61}{120},\
\frac{67}{120},\ \frac{71}{120},\ \frac{31}{120},\ \frac{101}{120},\
\frac{107}{120},\ \frac{59}{120},\ \frac{83}{120},\ \frac{43}{120},\
\frac{23}{120},\ \frac{103}{120}.
\end{equation}
After this, other gamma values at rational points with the denominator 120
are determined by reflection formulas.

Now we show the third statement of the Theorem. We apply Lemma \ref{kubert}
with $N=120$ and conclude that the generating set of 16 gamma values in
(\ref{gammabasis}) and (\ref{extravals}) cannot be made smaller. Therefore
the standard equations do not imply any relations between those 16 values.
This certainly applies to the subset (\ref{gammabasis}) of gamma values.
\qed

Explicit expressions of Chapter \ref{gammaevaluations} and of all gamma
values at rational points with the denominator dividing 120, in terms of
(\ref{gammabasis}) and (\ref{extravals}), are available as small {\sf
Mathematica} and {\sf Maple} packages via webpage \cite{homepage}. In those
packages, the list (\ref{algnumbers}) of algebraic numbers for
multiplicative expressions of the gamma values is extended by
\[
\sqrt6\pm\sqrt5,\quad \sqrt{10}\pm3,\quad \sqrt{\phi}\pm\sqrt5,\quad
\sqrt5\pm\sqrt{\phi^\star},\quad \sqrt{\phi}\pm\sqrt3,\quad
\sqrt3\pm\sqrt{\phi^\star}.
\]

Here we discuss briefly algebraic independence of gamma values. Rohrlich's
conjecture \cite[Section~3.3]{diophantine} implies that all multiplicative
relations between gamma values at rational points are implied by the
standard equations. The stronger Lang's conjecture \cite[Ch.~2]{lang}
implies that all algebraic relations between gamma values at rational points
are implied by the standard equations. If these conjectures are true, we
have (respectively) multiplicative or algebraic independence of the
distinguished gamma values in (\ref{gammabasis}). Analogues of these
conjectures are proved for the Thakur's gamma function in positive
characteristic \cite{gammap}. It is known that all $\overline{\QQ}$-linear
relations between the beta values ${\rm B}(a,b)$ with
$a,b,a+b\in\QQ\setminus\ZZ$ are implied by the standard equations
\cite{wolfart}. Chudnovsky proved that $\GAMMA13$ and $\GAMMA14$ are
transcendental over $\QQ(\pi)$; see \cite{waldschmidt}. Algebraicity of
gamma terms is determined, assuming Rohrlich's conjecture, by Koblitz-Ogus
criterion \cite{koblitzo}.

\section{Application to algebraic hypergeometric functions}
\label{application}

It is well-known that Gauss hypergeometric functions can be evaluated at
some points in terms of the gamma function. The most important formula is
the Gauss' identity \cite[Theorem 2.2.2]{specfaar}:
\begin{equation} \label{gaussid}
\hpg{2}{1}{a,\;b\;}{c}{\,1\,}=\frac{\Gamma(c)\,\Gamma(c-a-b)}{\Gamma(c-a)\,\Gamma(c-b)},
\qquad \mbox{if Re}(c-a-b)>0.
\end{equation}
Akin to this identity, linear relations between the 24 Kummer's solutions of
the hypergeometric differential equation (and analytic continuation formulas
for the Gauss hypergeometric function) involve similar gamma terms. An
example is \cite[Corollary 2.3.3]{specfaar}:
\begin{eqnarray*}
\hpg{2}{1}{a,\,b\,}{c}{\,x}&=&\frac{\Gamma(c)\,\Gamma(c-a-b)}{\Gamma(c-a)\,\Gamma(c-b)}
\;\hpg{2}{1}{a,\;b\,}{a+b+1-c}{\,1-x}\\
&&\!+\frac{\Gamma(c)\,\Gamma(a+b-c)}{\Gamma(a)\,\Gamma(b)}\,(1-x)^{c-a-b}
\;\hpg{2}{1}{c-a,\;c-b\,}{c+1-a-b}{\,1-x}.
\end{eqnarray*}

Our formulas in Section \ref{gammaevaluations} can be used to evaluate
widely used instances of the Gauss hypergeometric function most explicitly.
In particular, evaluation of algebraic Gauss hypergeometric functions with
(\ref{gaussid}) requires gamma values at rational points with the
denominator dividing 24 or 60; see \cite{schwartz}, \cite{kato},
\cite{alggauss}. The formulas in Section \ref{gammaevaluations} are
precisely sufficient to evaluate algebraic Gauss hypergeometric functions at
$x=1$. 
Here are a few examples:
\begin{eqnarray*}
\hpg{2}{1}{\frac14,\,-\frac1{12}\,}{\frac23}{\,1\,}\equal
\frac{\GAMMA12\GAMMA23}{\GAMMA34\GAMMA5{12}}\equall\frac{\sqrt{\sqrt3+1}}{2^{1/4}\,3^{3/8}},\\
\hpg{2}{1}{\frac5{24},\,-\frac1{24}\,}{\frac23}{\,1\,}\equal
\frac{\GAMMA12\GAMMA23}{\GAMMA{11}{24}\GAMMA{17}{24}}\equall\frac{\sqrt{\sqrt3+\sqrt2}}{2^{1/12}\,\sqrt3},\\
\hpg{2}{1}{\frac{11}{60},\,-\frac{1}{60}\,}{\frac23}{\,1\,}\equal
\frac{\GAMMA12\GAMMA23}{\GAMMA{29}{60}\GAMMA{41}{60}}\equall\frac{\sqrt{\phi^\star}\,\sqrt{\sqrt3+1}
\,\sqrt{\sqrt5+\sqrt3}}{2\,\sqrt3\;5^{7/24}},\\
\hpg{2}{1}{\frac3{10},\,-\frac1{30}\,}{\frac35}{\,1\,}\equal
\frac{\GAMMA13\GAMMA35}{\GAMMA{3}{10}\GAMMA{19}{30}}\equall
\frac{\sqrt{\sqrt{15}+\psi}}{2^{19/30}\,3^{1/20}\,5^{1/3}}.
\end{eqnarray*}
More of these evaluations are independently presented in \cite{kato}.
Explicit expressions for algebraic Gauss hypergeometric functions are given
in \cite{alggauss}.

\section{Gamma values and (hyper)elliptic integrals}
\label{gammaelliptic}

As is known \cite{borwzu}, some gamma values can be expressed in terms of
special values of the elliptic $\bf K$-function (i.e., {\em complete
elliptic integral of the first kind}). The advantage of doing this is that
numerical values of the elliptic $\bf K$-function can be very effectively
computed using the arithmetic-geometric mean; see
\cite[Section~3.2]{specfaar} and \cite[pg.~137]{borwbai}.

Here we present direct expressions for most of gamma values from
(\ref{gammabasis}) in terms of the elliptic $\bf K$-function. By combining
the data in \cite[pages \sf EllipticIntegralSingularValue.html,
EllipticLambdaFunction.html]{wolfram} and our formulas in Section
\ref{gammaevaluations} we get the following evaluations:
\begin{eqnarray*}
\GAMMA13\equal\frac{\pi^{1/3}\,2^{7/9}}{3^{1/12}}\;\KK{\frac{\sqrt{3}-1}{2\,\sqrt2}}^{1/3},\\
\GAMMA14\equal2\,\pi^{1/4}\;\KK{\frac{1}{\sqrt2}}^{1/2},\\
\GAMMA18\equal\pi^{1/8}\,2^{17/8}\;\KK{\frac{1}{\sqrt2}}^{1/4}\,\KK{\sqrt2-1}^{1/2},\\
\GAMMA1{15}\equal\frac{\pi^{1/6}\,3^{29/60}\,5^{1/24}\,\sqrt{\phi^\star}\,\sqrt{\psi+\sqrt3}}
{2^{1/9}}\;\GAMMA15^{1/2}\GAMMA25^{-1/2}\KK{\frac{\sqrt3-1}{2\,\sqrt2}}^{1/6}\\
&&\times\;
\KK{\frac{\left(2-\sqrt3\right)\left(3-\sqrt5\right)\left(\sqrt5-\sqrt3\right)}{8\,\sqrt2}}^{1/2},\\
\GAMMA1{20}\equal\frac{2^{9/40}\,5^{1/8}\,\phi^{5/8}\,\sqrt{\psi^\star+1}}{\pi^{1/4}}
\;\GAMMA15^{1/2}\GAMMA25^{1/2}\KK{\sqrt{\frac12-\sqrt{\sqrt5-2}}}^{1/2},\\
\GAMMA1{24}\equal\pi^{1/24}\,2^{89/36}\,3^{25/48}\,\sqrt{\sqrt2+1}\left(\sqrt3-1\right)^{1/4}
\,\KK{\frac{1}{\sqrt2}}^{1/4}\KK{\frac{\sqrt{3}-1}{2\,\sqrt2}}^{1/3}\\
&&\times\;\KK{(2-\sqrt3)(\sqrt3-\sqrt2)}^{1/2}.
\end{eqnarray*}

Other gamma values in (\ref{gammabasis}) can be expressed in terms of
hyperelliptic integrals. For example, consider the integrals:
\[
H_1=\int_0^1\frac{dz}{\sqrt{1-z^5}},\qquad\qquad
H_2=\int_0^1\frac{z\,dz}{\sqrt{1-z^5}}.
\]
Note that the two differentials under integration form a basis for the space
of holomorphic differentials on the genus 2 curve $y^2=x^5-1$. We have ${\rm
B}\left(\frac15,\frac12\right)=5\,H_1$, ${\rm
B}\left(\frac25,\frac12\right)=5\,H_2$; this can be seen after the
substitution $t\mapsto z^5$ in the standard definition
\cite[Definition~1.1.3]{specfaar} of the beta integral. From the two beta
values we conclude
\[
\GAMMA15=\pi^{1/5}\,2^{19/50}\,\sqrt5\,\phi^{1/10}\,H_1^{2/5}\,H_2^{1/5},\quad
\GAMMA25=\pi^{2/5}\,2^{4/25}\,\phi^{1/5}\,H_1^{-1/5}\,H_2^{2/5}.
\]

To get similar expressions for gamma values at rational points with the
denominator 60, one may substitute $t\mapsto z^{30}$ into the beta integrals
${\rm B}\left(\frac1{60},\,\frac12\right)$ and ${\rm
B}\left(\frac7{60},\,\frac12\right)$ and get hyperelliptic integrals on the
genus 15 curve $y^2=x^{31}-x$. Possibly, there are integrals on lower genus
curves related to those gamma values.

\end{document}